\documentclass[leqno,12pt]{amsart}
\usepackage{amsfonts}
\usepackage{amsmath,amssymb,amsthm}

\newcommand{\R}{\mathbb R}

\newcommand{\E}{\mathbb E}

\renewcommand{\span}{\mathrm{span}}
\newcommand{\tr}{\mathrm{tr}}

\newtheorem{thm}{Theorem}[section]

\newtheorem{prop}[thm]{Proposition}
\theoremstyle{definition}

\theoremstyle{remark}

\newcommand{\ds}{\displaystyle}

\begin{document}

\title[Quasi-minimal Rotational Surfaces]
{Quasi-minimal Rotational Surfaces in Pseudo-Euclidean
Four-dimensional Space}

\author{Georgi Ganchev and Velichka Milousheva}
\address{Bulgarian Academy of Sciences, Institute of Mathematics and Informatics,
Acad. G. Bonchev Str. bl. 8, 1113 Sofia, Bulgaria}
\email{ganchev@math.bas.bg}
\address{Bulgarian Academy of Sciences, Institute of Mathematics and Informatics,
Acad. G. Bonchev Str. bl. 8, 1113, Sofia, Bulgaria; "L. Karavelov"
Civil Engineering Higher School, 175 Suhodolska Str., 1373 Sofia,
Bulgaria} \email{vmil@math.bas.bg}


\begin{abstract}
In the four-dimensional pseudo-Euclidean space with neutral metric
there are three types of rotational surfaces with two-dimensional
axis -- rotational surfaces of elliptic, hyperbolic or parabolic type.
A surface whose mean curvature vector field is lightlike is said to be
quasi-minimal. In this paper we classify all rotational quasi-minimal
surfaces of elliptic, hyperbolic and parabolic type, respectively.
\end{abstract}

\maketitle

\section{Introduction}

A spacelike surface in the Minkowski 4-space $\E^4_1$ whose  mean curvature
vector $H$ is lightlike at each point is called \emph{marginally trapped}.
The concept of trapped surfaces,  introduced in 1965 by Roger Penrose
\cite{Pen}, plays an important role in general relativity and the theory
of cosmic black holes. Recently, classification results on marginally
trapped surfaces have been obtained imposing some extra conditions on
the mean curvature vector, the Gauss curvature or the second fundamental form.
In particular, marginally trapped surfaces with positive relative
nullity were classified by  B.-Y. Chen and J. Van der Veken in
\cite{Chen-Veken-1}. They also proved the non-existence of
marginally trapped surfaces in Robertson-Walker spaces with positive
relative nullity \cite{Chen-Veken-2} and classified marginally
trapped surfaces with parallel mean curvature vector in Lorenz space
forms \cite{Chen-Veken-3}. For a recent survey on marginally trapped
surfaces, see also \cite{Chen-book}.

In the four-dimensional Minkowski space there are three
types of rotational surfaces with two-dimensional axis -- rotational
surfaces of elliptic, hyperbolic or parabolic type, known also as
surfaces invariant under spacelike rotations, hyperbolic rotations
or screw rotations,  respectively.
A \emph{rotational surface of elliptic type} is an orbit of a regular
curve under the action of the orthogonal transformations of $\E^4_1$
which leave a timelike plane point-wise fixed. Similarly, a \emph{rotational
surface of hyperbolic type} is an orbit of a  regular curve under
the action of the orthogonal transformations of $\E^4_1$ which leave
a spacelike plane point-wise fixed. A \emph{rotational surface of
parabolic type} is an an orbit of a  regular curve under the action
of the orthogonal transformations of $\E^4_1$ which leave a
degenerate plane point-wise fixed.

The marginally trapped surfaces in Minkowski 4-space which are
invariant under spacelike rotations (rotational surfaces of elliptic
type) were classified by S. Haesen and  M. Ortega  in
\cite{Haesen-Ort-2}.   The classification of marginally trapped
surfaces   in $\R^4_1$ which are invariant under  boost
transformations (rotational surfaces of hyperbolic type) was
obtained  in \cite{Haesen-Ort-1} and the classification of
marginally trapped surfaces which are invariant under  screw
rotations (rotational surfaces of parabolic type) is given in
\cite{Haesen-Ort-3}.

Some classification results for rotational surfaces in three-dimensional
space forms satisfying some classical extra conditions have also been obtained.
For example, a classification of all timelike and
spacelike hyperbolic rotational surfaces with non-zero constant mean
curvature in the three-dimensional de Sitter space $\mathbb{S}^3_1$
is given in \cite{Liu-Liu} and a classification of the spacelike and
timelike Weingarten rotational surfaces of the three types in
$\mathbb{S}^3_1$ is found in \cite{Liu-Liu-2}. In \cite{GM3} we
described all Chen spacelike rotational surfaces of hyperbolic or
elliptic type.

Pseudo-Riemannian geometry has many important applications in
physics. According to the words of Bang-Yen Chen in his new book
\emph{Pseudo-Riemannian Geometry, $\delta$-Invariants and
Applications}, 2011: ''Spacetimes are the arenas in which all
physical events take place'' \cite{Chen-book}. In recent times,
physics and astrophysics have played a central role in shaping the
understanding of the universe through scientific observation and
experiment. The use of higher dimensional pseudo-Riemannian
manifolds in physics has led to many new developments in string
theory.

In the  pseudo-Riemannian geometry there is  an important subject
closely related with marginally trapped surfaces, namely
quasi-minimal surfaces. A Lorentz surface in a pseudo-Riemannian
manifold is called \emph{quasi-minimal}, if its mean curvature
vector is lightlike at each point of the surface. Borrowed from
general relativity, some authors call the quasi-minimal Lorentz
surfaces in a pseudo-Riemannian manifold also \emph{marginally
trapped}. We shall use the notion of a quasi-minimal surface.

The classification of quasi-minimal surfaces with parallel mean
curvature vector in the pseudo-Euclidean space $\E^4_2$ is obtained
in \cite{Chen-Garay}. In \cite{Chen1}  B.-Y. Chen classified
quasi-minimal Lorentz  flat surfaces in $\E^4_2$. As an application,
he gave the complete classification of biharmonic Lorentz surfaces
in  $\E^4_2$ with lightlike mean curvature vector. Several other
families of quasi-minimal surfaces  have also been classified. For
example, quasi-minimal surfaces with constant Gauss curvature in
$\E^4_2$ were classified in \cite{Chen2, Chen-Yang}. Quasi-minimal
Lagrangian surfaces and quasi-minimal slant surfaces in complex
space forms were classified, respectively, in \cite{Chen-Dillen} and
\cite{Chen-Mihai}. For an up-to-date survey on  quasi-minimal
surfaces, see also \cite{Chen3}.

In the present paper we consider three types of rotational surfaces
in the four-dimensional  pseudo-Euclidean space $\E^4_2$, namely
rotational surfaces of elliptic, hyperbolic, and parabolic type,
which are analogous to the three types of rotational surfaces in the
Minkowski space $\E^4_1$.

In Theorem \ref{T:quasi-minimal elliptic}  we find all
quasi-minimal rotational surfaces of elliptic type. In Theorem
\ref{T:quasi-minimal hyperbolic}  we describe all quasi-minimal
rotational surfaces of hyperbolic type and in Theorem
\ref{T:quasi-minimal parabolic}  we describe the construction of
all quasi-minimal rotational surfaces of parabolic type.

Our idea to study quasi-minimal rotational surfaces in the pseudo-Euclidean
space $\E^4_2$ was motivated by the results of S. Haesen and M. Ortega for
marginally trapped rotational surfaces in the Minkowski space
$\E^4_1$ \cite{Haesen-Ort-1, Haesen-Ort-2, Haesen-Ort-3}.

\section{Preliminaries}

Let  $\E^4_2$ be the pseudo-Euclidean 4-space endowed with the canonical pseudo-Euclidean metric of index 2 given by
$$g_0 = dx_1^2 + dx_2^2 - dx_3^2 - dx_4^2,$$
where $(x_1, x_2, x_3, x_4)$ is a rectangular coordinate system of $\E^4_2$. As usual, we denote by
$\langle \, , \rangle$ the indefinite inner scalar product with respect to $g_0$.

 A vector $v$ is called  \emph{spacelike} (respectively, \emph{timelike}) if $\langle v, v \rangle > 0$ (respectively, $\langle v, v \rangle < 0$).
 A vector $v$ is called \emph{lightlike} if it is nonzero and satisfies $\langle v, v \rangle = 0$.

A surface $M^2$ in $\E^4_2$ is called \emph{Lorentz}  if the
induced  metric $g$ on $M^2$ is Lorentzian. Thus at each point $p
\in M^2$ we have the following decomposition
$$\E^4_2 = T_pM^2 \oplus N_pM^2$$
with the property that the restriction of the metric
onto the tangent space $T_pM^2$ is of
signature $(1,1)$, and the restriction of the metric onto the normal space $N_pM^2$ is of signature $(1,1)$.

Denote by $\nabla$ and $\nabla'$ the Levi Civita connections of $M^2$  and $\E^4_2$, respectively.
Let $x$ and $y$ denote vector fields tangent to $M^2$ and let $\xi$ be a normal vector field.
The formulas of Gauss and Weingarten give a decomposition of the vector fields $\nabla'_xy$ and
$\nabla'_x \xi$ into a tangent and a normal component:
$$\begin{array}{l}
\vspace{2mm}
\nabla'_xy = \nabla_xy + \sigma(x,y);\\
\vspace{2mm}
\nabla'_x \xi = - A_{\xi} x + D_x \xi,
\end{array}$$
which define the second fundamental form $\sigma$, the normal
connection $D$, and the shape operator $A_{\xi}$ with respect to
$\xi$. In general, $A_{\xi}$ is not diagonalizable.

It is well known that the shape operator and the second fundamental form are related by the formula
$$\langle \sigma(x,y), \xi \rangle = \langle A_{\xi} x, y \rangle.$$

 The mean curvature vector  field $H$ of the surface $M^2$
is defined as $H = \ds{\frac{1}{2}\,  \tr\, \sigma}$.

A  surface $M^2$  is called \emph{minimal} if its mean curvature vector vanishes identically, i.e. $H =0$.
A natural extension of minimal surfaces are quasi-minimal surfaces.
The  surface $M^2$  is \emph{quasi-minimal} if its
mean curvature vector is lightlike at each point, i.e. $H \neq 0$ and $\langle H, H \rangle =0$.
Obviously, quasi-minimal surfaces are always non-minimal.

\section{Lorentz  rotational surfaces in pseudo-Euclidean 4-space}

 Let $Oe_1e_2e_3e_4$ be a fixed orthonormal coordinate system in the pseudo-Euclidean space $\E^4_2$, i.e. $\langle e_1, e_1 \rangle =
\langle e_2, e_2 \rangle  = 1, \, \langle e_3, e_3 \rangle =  \langle e_4, e_4 \rangle = -1$.

First we consider rotational surfaces of elliptic type.
Let $c: \widetilde{z} = \widetilde{z}(u), \,\, u \in J$ be a smooth spacelike curve,
parameterized by
$$\widetilde{z}(u) = \left( x_1(u), x_2(u),  r(u), 0 \right); \quad u \in J.$$
The curve $c$ lies in the three-dimensional subspace $\E^3_1 = \span\{e_1, e_2, e_3\}$ of $\E^4_2$.
Without loss of generality we assume that $c$ is
parameterized by the arc-length, i.e. $(x_1')^2 + (x_2')^2 - (r')^2 = 1$.
We assume also that  $r(u)>0, \,\, u \in J$.

Let us  consider the surface $\mathcal{M}'$ in $\E^4_2$ defined by
\begin{equation} \label{E:Eq-1}
\mathcal{M}': z(u,v) = \left( x_1(u), x_2(u), r(u) \cos v, r(u) \sin v\right);
\quad u \in J,\,\,  v \in [0; 2\pi).
\end{equation}

The tangent space of $\mathcal{M}'$ is spanned by the vector fields
$$\begin{array}{l}
\vspace{2mm}
z_u = \left(x_1', x_2', r' \cos v, r' \sin v  \right);\\
\vspace{2mm} z_v = \left( 0, 0 , - r \sin v, r \cos v\right).
\end{array}$$
Hence, the coefficients of the first fundamental form of $\mathcal{M}'$ are
$$E = \langle z_u, z_u \rangle = 1; \quad F = \langle z_u, z_v \rangle = 0; \quad G = \langle z_v, z_v \rangle = - r^2(u).$$

The surface $\mathcal{M}'$, defined by \eqref{E:Eq-1}, is a Lorentz surface  in $\E^4_2$, obtained by the rotation of the spacelike curve $c$ about
the two-dimensional Euclidean plane $Oe_1e_2$.
We  call $\mathcal{M}'$ a \emph{rotational surface of elliptic type}.

We can also obtain a rotational surface of elliptic type in  $\E^4_2$ using rotation of  a timelike curve about the two-dimensional plane $Oe_3e_4$.
Indeed, if $c$ is a timelike curve lying in the three-dimensional subspace $\span\{e_1, e_3, e_4\}$ of $\E^4_2$ and parameterized by
$$\widetilde{z}(u) = \left(r(u), 0, x_3(u), x_4(u) \right); \quad u \in J,$$
then the surface, defined by
$$z(u,v) = \left(r(u) \cos v, r(u) \sin v,  x_3(u), x_4(u) \right);
\quad u \in J,\,\,  v \in [0; 2\pi)$$
is a Lorentz rotational surface of elliptic type.

\vskip 2mm
Next, we consider rotational surfaces of hyperbolic type.
Let $c: \widetilde{z} = \widetilde{z}(u), \,\, u \in J$ be a smooth spacelike curve, lying in the three-dimensional subspace $\E^3_1 = \span\{e_1, e_2, e_4\}$ of $\E^4_2$
and parameterized by
$$\widetilde{z}(u) = \left( r(u), x_2(u), 0, x_4(u) \right); \quad u \in J.$$
Without loss of generality we assume that $c$ is
parameterized by the arc-length, i.e. $(r')^2 +(x_2')^2 - (x_4')^2  = 1$.
We assume also that  $r(u)>0, \,\, u \in J$.

Now we consider the surface $\mathcal{M}''$ in $\E^4_2$ defined by
\begin{equation} \label{E:Eq-2}
\mathcal{M}'': z(u,v) = \left(r(u) \cosh v,  x_2(u), r(u) \sinh v, x_4(u) \right);
\quad u \in J,\,\,  v \in \R.
\end{equation}

The tangent space of $\mathcal{M}''$ is spanned by the vector fields
$$\begin{array}{l}
\vspace{2mm}
z_u = \left( r' \cosh v, x_2',  r' \sinh v, x_4' \right);\\
\vspace{2mm} z_v = \left(r \sinh v, 0, r \cosh v, 0 \right),
\end{array}$$
and the coefficients of the first fundamental form of $\mathcal{M}''$ are
$$E = \langle z_u, z_u \rangle = 1; \quad F = \langle z_u, z_v \rangle = 0; \quad G = \langle z_v, z_v \rangle = - r^2(u).$$

The surface $\mathcal{M}''$, defined by \eqref{E:Eq-2}, is a Lorentz surface  in $\E^4_2$, obtained by hyperbolic rotation of the spacelike curve $c$ about
the two-dimensional Lorentz plane $Oe_2e_4$.
We  call $\mathcal{M}''$ a \emph{rotational surface of hyperbolic type}.

Similarly, we can obtain a rotational surface of hyperbolic type  using hyperbolic rotation of  a timelike curve lying in $\span\{e_2, e_3, e_4\}$
about the two-dimensional Lorentz plane $Oe_2e_4$.
Indeed, if $c$ is a timelike curve  parameterized by
$$\widetilde{z}(u) = \left( 0, x_2(u), r(u), x_4(u) \right); \quad u \in J,$$
then the surface, defined by
$$z(u,v) = \left( r(u) \sinh v, x_2(u), r(u) \cosh v, x_4(u) \right);
\quad u \in J,\,\,  v \in \R$$
is a Lorentz rotational surface of hyperbolic type.

Rotational surfaces of hyperbolic type can also be obtained by  hyperbolic rotations of spacelike or timelike curves about the two-dimensional Lorentz planes $Oe_1e_3$,
 $Oe_1e_4$ and $Oe_2e_3$. We are not going to define all of them here, since they are constructed in a similar way.

\vskip 2mm
Now, let us consider rotational surfaces of parabolic type in $\E^4_2$.
For convenience we shall use the pseudo-orthonormal base $\{e_1,
e_4, \xi_1, \xi_2 \}$  of $\E^4_2$, such that $ \ds{\xi_1=
\frac{e_2 + e_3}{\sqrt{2}}},\,\, \ds{\xi_2= \frac{ - e_2 +
e_3}{\sqrt{2}}}$. Note that
$$\langle \xi_1, \xi_1 \rangle =0; \quad \langle \xi_2, \xi_2 \rangle =0; \quad \langle \xi_1, \xi_2 \rangle = -1.$$

Let $c$ be a spacelike curve lying in the subspace $\E^3_1 = \span\{e_1, e_2, e_3\}$ of $\E^4_2$ and parameterized by
$$ \widetilde{z}(u) =  x_1(u)\, e_1 + x_2(u) \, e_2 + x_3(u) \, e_3; \quad u \in J,$$
or equivalently,
$$ \widetilde{z}(u) =  x_1(u)\, e_1 + \frac{x_2(u) + x_3(u)}{\sqrt{2}} \, \xi_1 + \frac{ - x_2(u) + x_3(u)}{\sqrt{2}} \, \xi_2; \quad u \in J.$$
Denote $f(u) = \ds{ \frac{x_2(u) + x_3(u)}{\sqrt{2}}}$, \; $g(u) = \ds{\frac{ - x_2(u) + x_3(u)}{\sqrt{2}}}$. Then
$$ \widetilde{z}(u) =  x_1(u)\, e_1 + f(u) \, \xi_1 + g(u) \, \xi_2.$$
Without loss of generality we assume that $c$ is
parameterized by the arc-length, i.e. $(x_1')^2 +(x_2')^2 - (x_3')^2  = 1$, or equivalently  $(x_1')^2 - 2f'g' = 1$.

We define a rotational surface of parabolic type
in the following way:
\begin{equation} \label{E:Eq-3}
\mathcal{M}''': z(u,v) = x_1(u)\, e_1 + f(u) \, \xi_1 + (-v^2 f(u) + g(u)) \, \xi_2 + \sqrt{2}\, v f(u) \, e_4;
\quad u \in J,\,\,  v \in \R.
\end{equation}

The tangent vector fields of $\mathcal{M}'''$ are
$$\begin{array}{l}
\vspace{2mm}
z_u = x_1'\, e_1 + \sqrt{2}\, v f'\, e_4 + f' \, \xi_1 + (-v^2 f' + g') \, \xi_2;\\
\vspace{2mm} z_v =\sqrt{2}\,  f\, e_4  - 2 v f \, \xi_2.
\end{array}$$
Hence, the coefficients of the first fundamental form of $\mathcal{M}'$ are
$$E = \langle z_u, z_u \rangle = 1; \quad F = \langle z_u, z_v \rangle = 0; \quad G = \langle z_v, z_v \rangle = - 2 f^2(u).$$

The surface $\mathcal{M}'''$, defined by \eqref{E:Eq-3}, is a Lorentz surface  in $\E^4_2$, which we call \emph{rotational surface of parabolic type}.
The rotational axis is the  two-dimensional plane spanned by $e_1$ (a spacelike vector field) and $\xi_1$ (a lightlike vector field).

\vskip 2mm
Similarly, we can obtain a rotational surface of parabolic type using a timelike curve lying in the subspace $\span\{e_2, e_3, e_4\}$ as follows.
Let $c$ be a timelike curve  given by
$$\widetilde{z}(u) = x_2(u)\, e_2 + x_3(u) \, e_3 + x_4(u) \, e_4; \quad u \in J.$$
We consider the lightlike vector fields $ \ds{\overline{\xi}_1=
\frac{e_2 + e_4}{\sqrt{2}}},\,\, \ds{\overline{\xi}_2= \frac{ - e_2 + e_4}{\sqrt{2}}}$. Then the parametrization of $c$ is expressed as
$$ \widetilde{z}(u) =  x_3(u)\, e_3 + \overline{f}(u) \, \overline{\xi}_1 + \overline{g}(u) \, \overline{\xi}_2,$$
where
 $\overline{f}(u) = \ds{ \frac{x_2(u) + x_4(u)}{\sqrt{2}}}$, \; $\overline{g}(u) = \ds{\frac{ - x_2(u) + x_4(u)}{\sqrt{2}}}$.

Now, let us  consider the surface defined as follows.
\begin{equation} \label{E:Eq-4}
z(u,v) =  \sqrt{2}\, v \overline{f}(u) \, e_1 +  x_3(u)\, e_3 + \overline{f}(u) \, \overline{\xi}_1 + (v^2 \overline{f}(u) + \overline{g}(u)) \, \overline{\xi}_2;
\quad u \in J,\,\,  v \in \R.
\end{equation}

The surface, given by \eqref{E:Eq-4}, is a Lorentz surface  in $\E^4_2$ whose coefficients of the first fundamental form are
$$E = \langle z_u, z_u \rangle = - 1; \quad F = \langle z_u, z_v \rangle = 0; \quad G = \langle z_v, z_v \rangle = 2 \overline{f}^2(u).$$

This surface is also a  \emph{rotational surface of parabolic type}, where the rotational axis is the  two-dimensional plane spanned by
$e_3$ (a timelike vector field) and $\overline{\xi}_1$ (a lightlike vector field).

\vskip 3mm
In what follows, we find all quasi-minimal surfaces in the three classes of rotational surfaces: elliptic type, hyperbolic type, and parabolic type.

\subsection{Quasi-minimal rotational surfaces of elliptic type}

Let us  consider the surface $\mathcal{M}'$ in $\E^4_2$ defined by \eqref{E:Eq-1}.
Since the generating curve $c$ is a spacelike curve
parameterized by the arc-length, i.e. $(x_1')^2 + (x_2')^2 - (r')^2 = 1$, then $(x_1')^2 + (x_2')^2 = 1 + (r')^2 $ and $x_1' x_1'' + x_2 ' x_2'' = r' r''$.
We shall use the following orthonormal tangent frame field:
$$X = z_u; \qquad Y = \ds{\frac{z_v}{r}},$$
and the normal frame field $\{n_1, n_2\}$, defined by
\begin{equation} \label{E:Eq-5}
\begin{array}{l}
\vspace{2mm}
n_1 = \displaystyle{\frac{1}{\sqrt{1+(r')^2}}\left(- x_2', x_1',0,0 \right)};\\
\vspace{2mm}
n_2 = \displaystyle{\frac{1}{\sqrt{1+(r')^2}} \left(r' x_1', r' x_2', (1+(r')^2) \cos v, (1+(r')^2) \sin v \right)}.
\end{array}
\end{equation}
Note that
$$\langle X, X \rangle = 1; \quad \langle X, Y \rangle = 0; \quad \langle Y, Y \rangle = -1;$$
$$\langle n_1, n_1 \rangle = 1; \quad \langle n_1, n_2 \rangle = 0; \quad \langle n _2, n_2 \rangle = -1.$$

The second
partial derivatives of $z(u,v)$ are expressed as follows
\begin{equation} \label{E:Eq-6}
\begin{array}{l}
\vspace{2mm}
z_{uu} = \left(x_1'', x_2'', r'' \cos v, r'' \sin v \right);\\
\vspace{2mm}
z_{uv} = \left(0, 0, - r' \sin v, r' \cos v \right);\\
\vspace{2mm} z_{vv} = \left(0, 0, - r \cos v, - r \sin v \right).
\end{array}
\end{equation}

By a straightforward computation, using \eqref{E:Eq-5} and  \eqref{E:Eq-6}, we obtain the components of the second fundamental form:
$$\begin{array}{ll}
\vspace{2mm}
\langle z_{uu}, n_1 \rangle = \displaystyle{\frac{1}{\sqrt{1+(r')^2}}(x_1' x_2'' - x_1'' x_2')}; & \qquad \langle z_{uu}, n_2 \rangle = \displaystyle{-\frac{r''}{\sqrt{1+(r')^2}}};\\
\vspace{2mm}
\langle z_{uv}, n_1 \rangle = 0; & \qquad  \langle z_{uv}, n_2 \rangle = 0;\\
\vspace{2mm}
\langle z_{vv}, n_1 \rangle = 0; & \qquad  \langle z_{vv}, n_2 \rangle = r \sqrt{1+(r')^2}.
\end{array}$$

Hence,
$$\begin{array}{l}
\vspace{2mm}
\sigma(z_u,z_u)=\ds{\frac{x_1' x_2'' - x_1'' x_2'}{\sqrt{1+(r')^2}} \, n_1 + \frac{r''}{\sqrt{1+(r')^2}}\, n_2},\\
\vspace{2mm}
\sigma(z_u,z_v)= 0,\\
\vspace{2mm}
\sigma(z_v,z_v) =\ds{\qquad \qquad \qquad \quad - r \sqrt{1+(r')^2}\, n_2}.
\end{array}$$
With respect to the orthonormal frame field $\{X,Y\}$ we get the formulas:
\begin{equation} \label{E:Eq-7}
\begin{array}{l}
\vspace{2mm}
\sigma(X,X)=\ds{\frac{x_1' x_2'' - x_1'' x_2'}{\sqrt{1+(r')^2}} \, n_1 + \frac{r''}{\sqrt{1+(r')^2}}\, n_2},\\
\vspace{2mm}
\sigma(X,Y)= 0,\\
\vspace{2mm}
\sigma(Y,Y) =\ds{\qquad \qquad \qquad \qquad - \frac{\sqrt{1+(r')^2}}{r}\, n_2}.
\end{array}
\end{equation}

Formulas \eqref{E:Eq-7} imply that the Gauss curvature $K$ of the rotational surface of elliptic  type $\mathcal{M}'$ is
\begin{equation} \label{E:Eq-7'}
K = \ds{- \frac{r''}{r}}
\end{equation}
and the normal mean curvature vector field $H$ is
\begin{equation}  \label{E:Eq-8}
H = \ds{\frac{1}{2 r\sqrt{1+(r')^2}} \left( r(x_1' x_2'' - x_1'' x_2') \, n_1 + (r r'' + (r')^2 + 1) \, n_2 \right)}.
\end{equation}

Equalities \eqref{E:Eq-7'} and \eqref{E:Eq-8} imply the following two statements.

\begin{prop}\label{P:flat elliptic}
The rotational surface of elliptic   type $\mathcal{M}'$ is flat if and only if $r'' = 0$.
\end{prop}

\begin{prop}\label{P:minimal elliptic}
The rotational surface of elliptic   type $\mathcal{M}'$  is minimal if and only if \\$x_1' x_2'' - x_1'' x_2' = 0$ \, and  \, $r r'' + (r')^2 + 1 = 0$.
\end{prop}

In the present paper we are interested in quasi-minimal rotational surfaces, so we assume that $(x_1' x_2'' - x_1'' x_2')^2 + (r r'' + (r')^2 + 1)^2 \neq 0$.

It follows from    \eqref{E:Eq-5} that
\begin{equation} \label{E:Eq-9}
\begin{array}{l}
\vspace{2mm}
\nabla'_X n_1 =\ds{- \frac{x_1' x_2'' - x_1'' x_2'}{\sqrt{1+(r')^2}} \, X + \frac{r'}{1+(r')^2}(x_1' x_2'' - x_1'' x_2') \, n_2},\\
\vspace{2mm}
\nabla'_Y n_1 = 0,\\
\vspace{2mm}
\nabla'_X n_2 =\ds{\frac{r''}{\sqrt{1+(r')^2}} \, X + \frac{r'}{1+(r')^2}(x_1' x_2'' - x_1'' x_2') \, n_1},\\
\vspace{2mm}
\nabla'_Y n_2 = \ds{\frac{\sqrt{1+(r')^2}}{r}\, Y}.
\end{array}
\end{equation}

\vskip 2mm
We can distinguish  two special classes of rotational surfaces of elliptic  type.

\vskip 2mm
I. Let  $x_1' x_2'' - x_1'' x_2' = 0$,\,  $r r'' + (r')^2 + 1 \neq 0$. In this case from the first two equalities of \eqref{E:Eq-9} we get
 $$\nabla'_X n_1 = 0;\quad \nabla'_Y n_1 = 0,$$
 which imply that the normal vector field $n_1$ is constant. Hence,
the rotational surface of elliptic  type $\mathcal{M}'$ lies in the hyperplane $\E^3_2$ of $\E^4_2$ orthogonal to $n_1$,
i.e. $\mathcal{M}'$ lies in the hyperplane $\E^3_2 =  \span \{X,Y,n_2\}$.

Moreover,  the mean curvature vector field of $\mathcal{M}'$ is:
\begin{equation} \notag
H = \ds{\frac{r r'' + (r')^2 + 1}{2 r\sqrt{1+(r')^2}} \, n_2}.
\end{equation}
Hence,  $\langle H, H \rangle = 0$ if and only if $H =0$ (i.e. $\mathcal{M}'$ is minimal).
Consequently, there are no quasi-minimal rotational surfaces of elliptic  type in the class
$x_1' x_2'' - x_1'' x_2' = 0$.

\vskip 2mm
II. Let $r r'' + (r')^2 + 1 = 0$, \, $x_1' x_2'' - x_1'' x_2' \neq 0$. In this case,  it can be proved that $\mathcal{M}'$ does not lie in any hyperplane of $\E^4_2$.
But, since the mean curvature vector field is
\begin{equation} \notag
H = \ds{\frac{x_1' x_2'' - x_1'' x_2'}{2 \sqrt{1+(r')^2}} \, n_1},
\end{equation}
we have again  that  $\langle H, H \rangle = 0$ if and only if $H =0$.
Consequently, there are no quasi-minimal rotational surfaces of elliptic  type in the class
$r r'' + (r')^2 + 1 = 0$.

\vskip 3mm
Further we shall consider general rotational surfaces of elliptic  type, i.e. we assume that   $x_1' x_2'' - x_1'' x_2' \neq 0$ and  $r r'' + (r')^2 + 1 \neq 0$
in an open interval $I \subset J$. In the next theorem we give a local description of all quasi-minimal rotational surfaces of elliptic  type.

\begin{thm}\label{T:quasi-minimal elliptic}
Given a smooth positive function $r(u): I \subset \R \rightarrow \R$, define the functions
$$\varphi(u) = \eta \int \ds{\frac{r r'' + (r')^2 + 1}{r (1+(r')^2)}} \, du, \quad \eta = \pm 1,$$
and
\begin{equation} \notag
\begin{array}{l}
\vspace{2mm}
x_1(u) = \int  \sqrt{1+(r')^2} \,\cos \varphi(u) \, du,\\
\vspace{2mm}
x_2(u) =\int  \sqrt{1+(r')^2} \,\sin \varphi(u) \, du.
\end{array}
\end{equation}
Then the spacelike curve $c: \widetilde{z}(u) = \left( x_1(u), x_2(u),  r(u), 0 \right)$ is a generating curve of a
quasi-minimal rotational surface of elliptic  type.

Conversely, any quasi-minimal rotational surface of elliptic  type is locally constructed as above.
\end{thm}

\noindent
\emph{Proof:}
Let $\mathcal{M}'$ be a general  rotational surface of elliptic  type generated by a spacelike curve
$c: \widetilde{z}(u) = \left( x_1(u), x_2(u),  r(u), 0 \right);\,\, u \in J$.
We assume that $c$ is parameterized by the arc-length  and $x_1' x_2'' - x_1'' x_2' \neq 0$, \, $r r'' + (r')^2 + 1 \neq 0$ for $u \in I \subset J$.

It follows from \eqref{E:Eq-8} that $\mathcal{M}'$ is quasi-minimal if and only if
$$r^2 (x_1' x_2'' - x_1'' x_2')^2 - (r r'' + (r')^2 + 1)^2 = 0,$$
or equivalently
\begin{equation} \label{E:Eq-11}
r (x_1' x_2'' - x_1'' x_2') = \eta (r r'' + (r')^2 + 1), \qquad \eta = \pm 1.
\end {equation}

Since the curve $c$  is parameterized by the arc-length, we have   $(x_1')^2 + (x_2')^2 = 1 + (r')^2$, which implies that
there exists a smooth function $\varphi = \varphi(u)$ such that
\begin{equation} \label{E:Eq-12}
\begin{array}{l}
\vspace{2mm}
x_1'(u) = \sqrt{1+(r')^2} \,\cos \varphi(u),\\
\vspace{2mm}
x_2'(u) = \sqrt{1+(r')^2} \,\sin \varphi(u).
\end{array}
\end{equation}
Using \eqref{E:Eq-12} we get $x_1' x_2'' - x_1'' x_2' = (1+(r')^2) \varphi'$. Hence, condition \eqref{E:Eq-11} for quasi-minimality of $\mathcal{M}'$
is written in terms of $r(u)$ and $\varphi(u)$ as follows:
\begin{equation} \label{E:Eq-13}
\varphi' = \eta \, \ds{\frac{r r'' + (r')^2 + 1}{r (1+(r')^2)}}.
\end{equation}

Consequently, the mean curvature vector field of a quasi-minimal  rotational surface of elliptic  type is given by the formula
$$H = \ds{\frac{r r'' + (r')^2 + 1}{2r \sqrt{1+(r')^2}}} (\eta\, n_1 + n_2).$$

Formula \eqref{E:Eq-13} allows us to recover $\varphi(u)$ from $r(u)$, up to integration constant.
Using formulas \eqref{E:Eq-12}, we can recover $x_1(u)$ and $x_2(u)$ from the functions $\varphi(u)$ and $r(u)$, up to integration constants.
Consequently, the quasi-minimal rotational surface of elliptic type $\mathcal{M}'$ is constructed as described in the theorem.

Conversely, if we are given a smooth function $r(u) > 0$, we can define the function
$$\varphi(u) = \eta \int \ds{\frac{r r'' + (r')^2 + 1}{r (1+(r')^2)}} \, du,$$
where $\eta = \pm 1$,  and consider the functions
\begin{equation} \notag
\begin{array}{l}
\vspace{2mm}
x_1(u) = \int  \sqrt{1+(r')^2} \,\cos \varphi(u) \, du,\\
\vspace{2mm}
x_2(u) =\int  \sqrt{1+(r')^2} \,\sin \varphi(u) \, du.
\end{array}
\end{equation}
A straightforward computation shows that the curve
$c: \widetilde{z}(u) = \left( x_1(u), x_2(u),  r(u), 0 \right)$ is a spacelike curve generating a quasi-minimal rotational surface of elliptic
type according to formula \eqref{E:Eq-1}.

\qed

\subsection{Quasi-minimal rotational surfaces of hyperbolic type}

Now, we shall consider the rotational surface of hyperbolic type  $\mathcal{M}''$ defined by \eqref{E:Eq-2}.
The generating curve $c$ is a spacelike curve
parameterized by the arc-length, i.e. $(r')^2 + (x_2')^2 - (x_4')^2  = 1$, and hence $(x_4')^2 - (x_2')^2 = (r')^2 - 1 $.
We assume that $(r')^2 \neq 1$, otherwise the surface lies in a 2-dimensional plane.
Denote by $\varepsilon$ the sign of $(r')^2 - 1$.

As in the elliptic case, we use the following orthonormal tangent frame field:
$$X = z_u; \qquad Y = \ds{\frac{z_v}{r}},$$
and the normal frame field $\{n_1, n_2\}$, defined by
\begin{equation} \label{E:Eq-5a}
\begin{array}{l}
\vspace{2mm}
n_1 = \displaystyle{\frac{1}{\sqrt{\varepsilon((r')^2 - 1)}}\left(0, x_4', 0, x_2' \right)};\\
\vspace{2mm}
n_2 = \displaystyle{\frac{1}{\sqrt{\varepsilon((r')^2 - 1)}} \left((1-(r')^2) \cosh v, - r' x_2',  (1-(r')^2) \sinh v, - r' x_4' \right)}.
\end{array}
\end{equation}
The orthonormal frame field $\{X, Y, n_1, n_2\}$ satisfies
$$\langle X, X \rangle = 1; \quad \langle X, Y \rangle = 0; \quad \langle Y, Y \rangle = -1;$$
$$\langle n_1, n_1 \rangle = \varepsilon; \quad \langle n_1, n_2 \rangle = 0; \quad \langle n _2, n_2 \rangle = - \varepsilon.$$

Calculating the second
partial derivatives of $z(u,v)$ we obtain
\begin{equation} \label{E:Eq-6a}
\begin{array}{l}
\vspace{2mm}
z_{uu} = \left(r'' \cosh v, x_2'',  r'' \sinh v,  x_4'' \right);\\
\vspace{2mm}
z_{uv} = \left(r' \sinh v, 0,  r' \cosh v, 0 \right);\\
\vspace{2mm} z_{vv} = \left(r \cosh v, 0,  r \sinh v, 0 \right).
\end{array}
\end{equation}

Formulas  \eqref{E:Eq-5a} and  \eqref{E:Eq-6a} imply  that the components of the second fundamental form of $\mathcal{M}''$ are:
$$\begin{array}{ll}
\vspace{2mm}
\langle z_{uu}, n_1 \rangle = \displaystyle{\frac{1}{\sqrt{\varepsilon((r')^2 - 1)}}(x_4' x_2'' - x_4'' x_2')}; & \qquad \langle z_{uu}, n_2 \rangle = \displaystyle{\frac{r''}{\sqrt{\varepsilon((r')^2 - 1)}}};\\
\vspace{2mm}
\langle z_{uv}, n_1 \rangle = 0; & \qquad  \langle z_{uv}, n_2 \rangle = 0;\\
\vspace{2mm}
\langle z_{vv}, n_1 \rangle = 0; & \qquad  \langle z_{vv}, n_2 \rangle = \ds{\frac{r(1- (r')^2)}{\sqrt{\varepsilon((r')^2 - 1)}}}.
\end{array}$$

Hence,  we obtain  the following formulas for the second fundamental form $\sigma$:
\begin{equation} \label{E:Eq-7a}
\begin{array}{l}
\vspace{2mm}
\sigma(X,X)=\ds{\frac{\varepsilon (x_4' x_2'' - x_4'' x_2')}{\sqrt{\varepsilon((r')^2 - 1)}} \, n_1 - \frac{\varepsilon r''}{\sqrt{\varepsilon((r')^2 - 1)}}\, n_2},\\
\vspace{2mm}
\sigma(X,Y)= 0,\\
\vspace{2mm}
\sigma(Y,Y) =\ds{\qquad \qquad \qquad \qquad + \frac{\varepsilon ((r')^2 - 1)}{r \sqrt{\varepsilon((r')^2 - 1)}}\, n_2}.
\end{array}
\end{equation}

Formulas \eqref{E:Eq-7a} imply that the Gauss curvature $K$ of the rotational surface of hyperbolic  type $\mathcal{M}''$ is
\begin{equation} \label{E:Eq-7'a}
K = \ds{- \frac{r''}{r}}
\end{equation}
and the normal mean curvature vector field $H$ is
\begin{equation}  \label{E:Eq-8a}
H = \ds{\frac{\varepsilon}{2 r\sqrt{\varepsilon((r')^2 - 1)}} \left( r(x_4' x_2'' - x_4'' x_2') \, n_1 - (r r'' + (r')^2 - 1) \, n_2 \right)}.
\end{equation}

The next two statements follow directly from equalities \eqref{E:Eq-7'a} and \eqref{E:Eq-8a}.

\begin{prop}\label{P:flat hyperbolic}
The rotational surface of hyperbolic  type $\mathcal{M}''$ is flat if and only if $r'' = 0$.
\end{prop}

\begin{prop}\label{P:minimal hyperbolic}
The rotational surface of hyperbolic   type $\mathcal{M}''$  is minimal if and only if $x_4' x_2'' - x_4'' x_2' = 0$ \, and \, $r r'' + (r')^2 - 1 = 0$.
\end{prop}

We assume that $(x_4' x_2'' - x_4'' x_2')^2 + (r r'' + (r')^2 -1)^2 \neq 0$, since we are interested in quasi-minimal rotational surfaces.

Similarly to the elliptic case it follows from    \eqref{E:Eq-5a} that
\begin{equation} \label{E:Eq-9a}
\begin{array}{l}
\vspace{2mm}
\nabla'_X n_1 =\ds{\frac{x_2' x_4'' - x_2'' x_4'}{\sqrt{\varepsilon ((r')^2 - 1)}} \, X + \frac{r'}{\varepsilon ((r')^2 - 1)}(x_2' x_4'' - x_2'' x_4') \, n_2},\\
\vspace{2mm}
\nabla'_Y n_1 = 0.
\end{array}
\end{equation}

\vskip 2mm
We  distinguish the following  two special classes of rotational surfaces of hyperbolic  type.

\vskip 2mm
I. Let  $x_2' x_4'' - x_2'' x_4' = 0$,\,  $r r'' + (r')^2 -1 \neq 0$. Using \eqref{E:Eq-9a} we get that in this case
 $$\nabla'_X n_1 = 0;\quad \nabla'_Y n_1 = 0,$$
 which imply that the rotational surface of hyperbolic  type $\mathcal{M}''$ lies in the hyperplane $\span \{X,Y,n_2\}$.

The mean curvature vector field of $\mathcal{M}''$ is:
\begin{equation} \notag
H = \ds{\frac{\varepsilon (1 - (r')^2 - r r'')}{2 r\sqrt{\varepsilon ((r')^2 - 1)}} \, n_2}.
\end{equation}
Hence,  $\langle H, H \rangle = 0$ if and only if $H =0$.
Consequently, there are no quasi-minimal rotational surfaces of hyperbolic  type in the class
$x_2' x_4'' - x_2'' x_4' = 0$.

\vskip 2mm
II. Let $r r'' + (r')^2 - 1 = 0$, \, $x_2' x_4'' - x_2'' x_4' \neq 0$.
In this case  $\mathcal{M}''$ does not lie in any hyperplane of $\E^4_2$ and
the mean curvature vector field is
\begin{equation} \notag
H = \ds{\frac{\varepsilon (x_4' x_2'' - x_4'' x_2')}{2 \sqrt{\varepsilon ((r')^2 - 1)}} \, n_1}.
\end{equation}
Hence, we have again  that  $\langle H, H \rangle = 0$ if and only if $H =0$.
Consequently, there are no quasi-minimal rotational surfaces of hyperbolic  type in the class
$r r'' + (r')^2 - 1 = 0$.

\vskip 3mm
Further we consider general rotational surfaces of hyperbolic  type, i.e. we assume that   $x_2' x_4'' - x_2'' x_4' \neq 0$ and  $r r'' + (r')^2 - 1 \neq 0$
in an open interval $I \subset J$. The following theorem gives a local description of all quasi-minimal rotational surfaces of hyperbolic  type.

\begin{thm}\label{T:quasi-minimal hyperbolic}
\emph{Case (A)}. Given a smooth positive function $r(u): I \subset \R \rightarrow \R$, such that $(r')^2 > 1$,  define the functions
$$\varphi(u) = \eta \int \ds{\frac{r r'' + (r')^2 - 1}{r (1-(r')^2)}} \, du, \quad \eta = \pm 1,$$
and
\begin{equation} \notag
\begin{array}{l}
\vspace{2mm}
x_2(u) = \int  \sqrt{(r')^2-1} \,\sinh \varphi(u) \, du,\\
\vspace{2mm}
x_4(u) =\int  \sqrt{(r')^2-1} \,\cosh \varphi(u) \, du.
\end{array}
\end{equation}
Then the spacelike curve $c: \widetilde{z}(u) = \left( r(u), x_2(u), 0, x_4(u) \right)$ is a generating curve of a
quasi-minimal rotational surface of hyperbolic type.

\emph{Case (B)}. Given a smooth positive function $r(u): I \subset \R \rightarrow \R$,  such that $(r')^2 < 1$, define the functions
$$\varphi(u) = \eta \int \ds{\frac{r r'' + (r')^2 - 1}{r (1-(r')^2)}} \, du, \quad \eta = \pm 1,$$
and
\begin{equation} \notag
\begin{array}{l}
\vspace{2mm}
x_2(u) = \int  \sqrt{1 - (r')^2} \,\cosh \varphi(u) \, du,\\
\vspace{2mm}
x_4(u) =\int  \sqrt{1 - (r')^2} \,\sinh \varphi(u) \, du.
\end{array}
\end{equation}
Then the spacelike curve $c: \widetilde{z}(u) = \left( r(u), x_2(u), 0, x_4(u) \right)$ is a generating curve of a
quasi-minimal rotational surface of hyperbolic type.

Conversely, any quasi-minimal rotational surface of hyperbolic type is locally described by one of the cases given above.
\end{thm}

\noindent
\emph{Proof:}
Let $\mathcal{M}''$ be a general  rotational surface of hyperbolic  type generated by a spacelike curve
$c: \widetilde{z}(u) = \left( r(u), x_2(u), 0, x_4(u) \right);\,\, u \in J$.
We assume that $c$ is parameterized by the arc-length  and $x_2' x_4'' - x_2'' x_4' \neq 0$, \, $r r'' + (r')^2 - 1 \neq 0$  in an interval $I \subset J$.

Formula \eqref{E:Eq-8a} implies that $\mathcal{M}''$ is quasi-minimal if and only if
\begin{equation} \label{E:Eq-11a}
r (x_2' x_4'' - x_2'' x_4') = \eta (r r'' + (r')^2 - 1), \qquad \eta = \pm 1.
\end {equation}

Since our considerations are local, we can assume that either $\varepsilon = 1$ in some open interval $I_0 \subset I$ or $\varepsilon = - 1$  in an open interval $I_1 \subset I$.
We study the restriction of $\mathcal{M}''$ on $I_0$, respectively $I_1$.

If $\varepsilon = 1$, then using that $(x_4')^2 - (x_2')^2 = (r')^2 - 1$ and $(r')^2 - 1 >0$ we obtain that
there exists a smooth function $\varphi = \varphi(u)$ such that
\begin{equation} \label{E:Eq-12a}
\begin{array}{l}
\vspace{2mm}
x_2'(u) = \sqrt{(r')^2 - 1} \,\sinh \varphi(u),\\
\vspace{2mm}
x_4'(u) = \sqrt{(r')^2 - 1} \,\cosh \varphi(u).
\end{array}
\end{equation}
The last equalities imply  $x_4' x_2'' - x_4'' x_2' = ((r')^2 - 1) \varphi'$. Hence, condition \eqref{E:Eq-11a} for quasi-minimality of $\mathcal{M}''$
is written in terms of $r(u)$ and $\varphi(u)$ as follows:
\begin{equation} \label{E:Eq-13a}
\varphi' = \eta \, \ds{\frac{1 - (r')^2 - r r''}{r ((r')^2 - 1)}}.
\end{equation}

Then, the mean curvature vector field  is given by the formula
$$H = \ds{\frac{1 - (r')^2 - r r''}{2r \sqrt{(r')^2-1}}} (\eta\, n_1 + n_2).$$

Using \eqref{E:Eq-13a} we can recover $\varphi(u)$ from $r(u)$, up to integration constant, and using  \eqref{E:Eq-12a},
we can recover $x_2(u)$ and $x_4(u)$ from the functions $\varphi(u)$ and $r(u)$, up to integration constants.
Consequently, if $\varepsilon = 1$  the restriction of the quasi-minimal rotational surface of hyperbolic type $\mathcal{M}''$  on  $I_0$
is constructed as described in case (A) of the theorem.

If $\varepsilon = - 1$, then
there exists a smooth function $\varphi = \varphi(u)$ such that
\begin{equation} \notag
\begin{array}{l}
\vspace{2mm}
x_2'(u) = \sqrt{1 - (r')^2} \,\cosh \varphi(u),\\
\vspace{2mm}
x_4'(u) = \sqrt{1 - (r')^2} \,\sinh \varphi(u).
\end{array}
\end{equation}
As in the previous case we get that condition \eqref{E:Eq-11a} for quasi-minimality of $\mathcal{M}''$
is:
\begin{equation} \notag
\varphi' = \eta \, \ds{\frac{r r'' + (r')^2 - 1}{r (1 - (r')^2)}},
\end{equation}
and the mean curvature vector field  is given by the formula
$$H = \ds{\frac{r r'' + (r')^2 - 1}{2r \sqrt{1-(r')^2}}} (\eta\, n_1 + n_2).$$

Hence, we can recover $\varphi(u)$ from $r(u)$, and
$x_2(u)$,  $x_4(u)$ from  $\varphi(u)$ and $r(u)$, up to integration constants.
Consequently, if $\varepsilon = - 1$  the restriction of the quasi-minimal rotational surface of hyperbolic type $\mathcal{M}''$  on  $I_1$
is constructed as described in case (B) of the theorem.

\vskip 2mm
Conversely, if we are given a smooth function $r(u) > 0$, we can define the function
$$\varphi(u) = \eta \int \ds{\frac{r r'' + (r')^2 - 1}{r (1-(r')^2)}} \, du, \quad \eta = \pm 1,$$
and consider the functions
$$\begin{array}{l}
\vspace{2mm}
x_2(u) = \int  \sqrt{(r')^2-1} \,\sinh \varphi(u) \, du,\\
\vspace{2mm}
x_4(u) =\int  \sqrt{(r')^2-1} \,\cosh \varphi(u) \, du,
\end{array}  \eqno{\textrm{case (A)}}$$
or
$$\begin{array}{l}
\vspace{2mm}
x_2(u) = \int  \sqrt{1 - (r')^2} \,\cosh \varphi(u) \, du,\\
\vspace{2mm}
x_4(u) =\int  \sqrt{1 - (r')^2} \,\sinh \varphi(u) \, du.
\end{array}  \eqno{\textrm{case (B)}}$$

A straightforward computation shows that the curve $c: \widetilde{z}(u) = \left( r(u), x_2(u), 0, x_4(u) \right)$
is a spacelike curve generating a quasi-minimal rotational surface of  hyperbolic type according to formula \eqref{E:Eq-2}.

\qed

\subsection{Quasi-minimal rotational surfaces of parabolic type}

Now we shall  consider the rotational surface of parabolic type $\mathcal{M}'''$ in $\E^4_2$ defined by formula \eqref{E:Eq-3} with respect to
 $\{e_1, e_4, \xi_1, \xi_2 \}$, where $ \ds{\xi_1= \frac{e_2 + e_3}{\sqrt{2}}},\,\, \ds{\xi_2= \frac{ - e_2 + e_3}{\sqrt{2}}}$. Recall that
$$\langle \xi_1, \xi_1 \rangle =0; \quad \langle \xi_2, \xi_2 \rangle =0; \quad \langle \xi_1, \xi_2 \rangle = -1.$$
The generating curve $c$ is a spacelike curve
parameterized by the arc-length, i.e. $(x_1')^2 + (x_2')^2 - (x_3')^2 = 1$, and hence $(x_1')^2 = 1 + 2 f' g' $; \, $x_1' x_1'' = g' f'' + f' g''$.

We  use the following orthonormal tangent frame field:
$$\begin{array}{l}
\vspace{2mm}
X = z_u = x_1'\, e_1 + \sqrt{2}\, v f'\, e_4 + f' \, \xi_1 + (-v^2 f' + g') \, \xi_2;\\
\vspace{2mm}
Y = \ds{\frac{z_v}{\sqrt{2}f} = e_4  - \sqrt{2} v \, \xi_2};
\end{array}$$
and the normal frame field $\{n_1, n_2\}$, defined by
\begin{equation} \label{E:Eq-5b}
\begin{array}{l}
\vspace{2mm}
n_1 = \displaystyle{e_1 + \frac{x_1'}{f'}\, \xi_2};\\
\vspace{2mm}
n_2 = \displaystyle{x_1'\, e_1 + \sqrt{2}\, v f'\, e_4 + f' \, \xi_1 + \frac{1 + f'g' - v^2 (f')^2}{f'} \, \xi_2}.
\end{array}
\end{equation}

The second
partial derivatives of $z(u,v)$ are expressed as follows
\begin{equation} \label{E:Eq-6b}
\begin{array}{l}
\vspace{2mm}
z_{uu} = x_1''\, e_1 + \sqrt{2}\, v f''\, e_4 + f'' \, \xi_1 + (-v^2 f'' + g'')  \, \xi_2;\\
\vspace{2mm}
z_{uv} = \sqrt{2} f' \,e_4  - 2 v f' \, \xi_2;\\
\vspace{2mm} z_{vv} = - 2 f \, \xi_2.
\end{array}
\end{equation}

By a straightforward computation from \eqref{E:Eq-5b} and  \eqref{E:Eq-6b} we obtain the components of the second fundamental form:
$$\begin{array}{ll}
\vspace{2mm}
\langle z_{uu}, n_1 \rangle = \displaystyle{\frac{x_1'' f' - x_1' f''}{f'}}; & \qquad \langle z_{uu}, n_2 \rangle = \displaystyle{-\frac{f''}{f'}};\\
\vspace{2mm}
\langle z_{uv}, n_1 \rangle = 0; & \qquad  \langle z_{uv}, n_2 \rangle = 0;\\
\vspace{2mm}
\langle z_{vv}, n_1 \rangle = 0; & \qquad  \langle z_{vv}, n_2 \rangle = 2f f'.
\end{array}$$

Hence, we obtain the following formulas for the second fundamental form $\sigma$:

\begin{equation} \label{E:Eq-7b}
\begin{array}{l}
\vspace{2mm}
\sigma(X,X)=\ds{\frac{x_1'' f' - x_1' f''}{f'} \, n_1 + \frac{f''}{f'}\, n_2},\\
\vspace{2mm}
\sigma(X,Y)= 0,\\
\vspace{2mm}
\sigma(Y,Y) =\ds{\qquad \qquad \qquad \quad - \frac{f'}{f}\, n_2}.
\end{array}
\end{equation}

Formulas \eqref{E:Eq-7b} imply that the Gauss curvature $K$ of the rotational surface of parabolic type $\mathcal{M}'''$ is expressed as
\begin{equation} \label{E:Eq-7'b}
K = \ds{- \frac{f''}{f}}
\end{equation}
and the mean curvature vector field $H$ is
\begin{equation}  \label{E:Eq-8b}
H = \ds{\frac{1}{2 f f'} \left( f(x_1'' f' - x_1' f'') \, n_1 + (f f'' + (f')^2) \, n_2 \right)}.
\end{equation}

Using equalities  \eqref{E:Eq-7'b} and \eqref{E:Eq-8b} we get the following two statements.

\begin{prop}\label{P:flat parabolic}
The rotational surface of parabolic   type $\mathcal{M}'''$ is flat if and only if $f'' = 0$.
\end{prop}

\begin{prop}\label{P:minimal parabolic}
The rotational surface of  parabolic   type $\mathcal{M}'''$  is minimal if and only if $x_1'' f' - x_1' f'' = 0$ \, and \, $f f'' + (f')^2 = 0$.
\end{prop}

We assume that $(x_1'' f' - x_1' f'')^2 + (f f'' + (f')^2)^2 \neq 0$, since  we study quasi-minimal rotational surfaces.

It follows from    \eqref{E:Eq-5b} that
\begin{equation} \label{E:Eq-9b}
\begin{array}{l}
\vspace{2mm}
\nabla'_X n_1 =\ds{- \frac{x_1'' f' - x_1' f''}{f'} \, X + \frac{x_1'' f' - x_1' f''}{f'} \, n_2},\\
\vspace{2mm}
\nabla'_Y n_1 = 0,\\
\vspace{2mm}
\nabla'_X n_2 =\ds{\frac{f''}{f'} \, X + \frac{x_1'' f' - x_1' f''}{f'}\, n_1},\\
\vspace{2mm}
\nabla'_Y n_2 = \ds{\frac{f'}{f}\, Y}.
\end{array}
\end{equation}

\vskip 2mm
As in the elliptic and hyperbolic cases we distinguish  two special classes of rotational surfaces of parabolic  type.

\vskip 2mm
I. Let  $x_1'' f' - x_1' f'' = 0$,\,  $f f'' + (f')^2 \neq 0$. In this case from the first two equalities of \eqref{E:Eq-9b} we get
 $$\nabla'_X n_1 = 0;\quad \nabla'_Y n_1 = 0,$$
 which imply that the normal vector field $n_1$ is constant and hence,
the rotational surface of parabolic  type $\mathcal{M}'''$ lies in the hyperplane $\E^3_2 =  \span \{X,Y,n_2\}$ of $\E^4_2$.

In this case  the mean curvature vector field of $\mathcal{M}'''$ is:
\begin{equation} \notag
H = \ds{\frac{f f'' + (f')^2}{2 f f'} \, n_2},
\end{equation}
which implies that $\langle H, H \rangle = 0$ if and only if $H =0$.
Consequently, there are no quasi-minimal rotational surfaces of parabolic type in this  class.

\vskip 2mm
II. Let $f f'' + (f')^2 = 0$, \, $x_1'' f' - x_1' f'' \neq 0$. In this case
the mean curvature vector field is
\begin{equation} \notag
H = \ds{\frac{x_1'' f' - x_1' f''}{2 f'} \, n_1},
\end{equation}
which implies again  that  $\langle H, H \rangle = 0$ if and only if $H =0$.
Consequently, there are no quasi-minimal rotational surfaces of parabolic  type in this special class.

\vskip 3mm
Further we consider general rotational surfaces of parabolic type, i.e. we assume that   $x_1'' f' - x_1' f'' \neq 0$ and  $f f'' + (f')^2 \neq 0$
in an open interval $I \subset J$. In the following theorem we give a local description of all quasi-minimal rotational surfaces of parabolic  type.

\begin{thm}\label{T:quasi-minimal parabolic}
Given a smooth function $f(u): I \subset \R \rightarrow \R$,
define the functions
$$\varphi(u) = \ds{f'(u) \left(C + \eta \left( - \frac{1}{f'(u)} + \int\frac{du}{f(u)} \right) \right)}, \quad \eta = \pm 1, \,\, C = const,$$
and
$$x_1(u) = \ds{\int\varphi(u) du}; \qquad g(u) = \ds{\int\frac{\varphi^2(u) - 1}{2 f'(u)} \,du}.$$
Then the curve $c: \widetilde{z}(u) =  x_1(u)\, e_1 + f(u) \, \xi_1 + g(u) \, \xi_2$ is a spacelike curve generating a
quasi-minimal rotational surface of parabolic  type.

Conversely, any quasi-minimal rotational surface of parabolic type is locally constructed as described above.
\end{thm}

\noindent
\emph{Proof:}
Let $\mathcal{M}'''$ be a general  rotational surface of parabolic  type generated by a spacelike curve
$c: \widetilde{z}(u) =  x_1(u)\, e_1 + f(u) \, \xi_1 + g(u) \, \xi_2; \,\, u \in J$.
We assume that $c$ is parameterized by the arc-length  and $x_1'' f' - x_1' f'' \neq 0$, \,  $f f'' + (f')^2 \neq 0$ for $u \in I \subset J$.

Equality \eqref{E:Eq-8b} implies that $\mathcal{M}'''$ is quasi-minimal if and only if
\begin{equation} \label{E:Eq-11b}
f (x_1'' f' - x_1' f'') = \eta (f f'' + (f')^2), \qquad \eta = \pm 1.
\end {equation}

Hence, the mean curvature vector field of a quasi-minimal  rotational surface of parabolic type is given by the formula
$$H = \ds{\frac{1}{2} \left( \ln |f f'| \right)'(\eta\, n_1 + n_2)}.$$

We denote $\varphi(u) = x_1'(u)$. Since $c$  is parameterized by the arc-length, we have   $(x_1')^2 = 1 + 2 f' g'$, which implies that
$g'(u) = \ds{\frac{\varphi^2(u)- 1}{2 f'(u)}}$. The last equality allows us to recover $g(u)$ from the functions $\varphi(u)$ and $f(u)$, up to integration constant.

Condition \eqref{E:Eq-11b} for quasi-minimality of $\mathcal{M}'''$
is written in terms of $f(u)$ and $\varphi(u)$ as follows:
\begin{equation} \label{E:Eq-13b}
\varphi'  - \frac{f''}{f'} \varphi = \eta \left(\frac{f''}{f'} + \frac{f'}{f} \right).
\end{equation}

We consider \eqref{E:Eq-13b} as a differential equation  with respect to $\varphi(u)$. Then the general solution of \eqref{E:Eq-13b} is given by the formula
\begin{equation} \label{E:Eq-31}
\varphi(u) = \ds{e^{- \int p(u) du}\left(C + \int q(u) \, e^{\int p(u) du} du \right)},
\end{equation}
where $p(u) = \ds{- \frac{f''}{f'}}$, \,\, $q(u) = \ds{ \eta \left(\frac{f''}{f'} + \frac{f'}{f} \right)}$.
Calculating the integrals in formula \eqref{E:Eq-31} we obtain
\begin{equation} \label{E:Eq-32}
\varphi(u) = \ds{f'(u) \left(C + \eta \left( - \frac{1}{f'(u)} + \int\frac{du}{f(u)} \right) \right)}, \quad \eta = \pm 1, \,\, C = const,
\end{equation}
which allows us to recover $\varphi(u)$ from $f(u)$.

Hence, the quasi-minimal rotational surface of parabolic type $\mathcal{M}'''$ is locally constructed as described in the theorem.

Conversely, if we are given a smooth function $f(u)$, we can define the function $\varphi(u)$ by formula \eqref{E:Eq-32} and
 consider the functions
\begin{equation} \notag
x_1 (u) = \ds{\int \varphi(u)\, du}; \qquad    g(u) = \int \ds{\frac{\varphi^2(u)- 1}{2 f'(u)} \, du}.
\end{equation}
A straightforward computation shows that the curve
$c: \widetilde{z}(u) =  x_1(u)\, e_1 + f(u) \, \xi_1 + g(u) \, \xi_2$ is a spacelike curve generating a
quasi-minimal rotational surface of parabolic type according to formula \eqref{E:Eq-3}.

\qed

\end{document}